\documentclass{amsart}
\usepackage{amssymb, amscd, amsmath, amsfonts}
\usepackage[dvipsnames]{xcolor}
\usepackage{verbatim}
\usepackage{geometry}                
\usepackage{graphicx}
\usepackage{epstopdf}
\newtheorem{theorem}{Theorem}
\newtheorem{prop}[theorem]{Proposition}

\newtheorem{definition}[theorem]{Definition}
\newtheorem{example}[theorem]{Example}

\def\la{{\lambda}}
\def\al{{\alpha}}
\def\ga{{\gamma}}

\def\QQ{{\mathbb Q}}

\def\st{{\tilde s}}

\def\hht{{\tilde h}}

\def\mt{{\tilde m}}

\def\bfp{{\mathbf p}}

\def\mvdash{{\,\vdash\!\!\vdash}}

\newdimen\squaresize \squaresize=10pt
\newdimen\thickness \thickness=0.4pt

\def\square#1{\hbox{\vrule width \thickness
     \vbox to \squaresize{\hrule height \thickness\vss
        \hbox to \squaresize{\hss#1\hss}
     \vss\hrule height\thickness}
\unskip\vrule width \thickness}
\kern-\thickness}

\def\vsquare#1{\vbox{\square{$#1$}}\kern-\thickness}

\def\young#1{
\vbox{\smallskip\offinterlineskip
\halign{&\vsquare{##}\cr #1}}}

\def\thisbox#1{\kern-.09ex\fbox{#1}}
\def\downbox#1{\lower1.200em\hbox{#1}}
\newdimen\Squaresize \Squaresize=15pt
\newdimen\Thickness \Thickness=0.4pt

\def\Square#1{\hbox{\vrule width \Thickness
     \vbox to \Squaresize{\hrule height \Thickness\vss
        \hbox to \Squaresize{\hss#1\hss}
     \vss\hrule height\Thickness}
\unskip\vrule width \Thickness}
\kern-\Thickness}

\def\Vsquare#1{\vbox{\Square{$#1$}}\kern-\Thickness}

\address{Dartmouth College, Hanover, NH, United States\\
and York University, Toronto, ON, Canada
}
\begin{document}

\title[Character symmetric functions]{Symmetric group characters
as symmetric functions (Extended abstract)}
\author[R. Orellana \and M. Zabrocki]{Rosa Orellana \and Mike Zabrocki}
\date{\today}

\maketitle

\begin{abstract}
The irreducible characters of the symmetric group
are a symmetric polynomial
in the eigenvalues of a permutation matrix.
They can therefore be realized as a symmetric function
that can be evaluated at a set of variables and form
a basis of the symmetric functions.  This basis of the
symmetric functions is of non-homogeneous degree and the
(outer) product structure coefficients are the
stable Kronecker coefficients.

We introduce the irreducible character basis by
defining it in terms of the induced
trivial characters of the symmetric group which also 
form a basis of the symmetric functions.  The irreducible
character basis is closely related to character polynomials
and we obtain some of the change of basis coefficients
by making this connection explicit.  Other change of
basis coefficients come from a representation theoretic
connection with the partition algebra,
and still others are derived by developing combinatorial
expressions.

This document is an extended abstract which can be used
as a review reference so that this basis can be implemented in Sage.
A more full version of the results in this abstract can be
found in {\tt arXiv:1605.06672}.
\end{abstract}


\begin{section}{Introduction}
We begin with a very basic question in representation theory.
Can the irreducible characters of the symmetric group be
realized as functions of the eigenvalues of the permutation
matrices?
In general, characters of matrix groups are
symmetric functions in the set of
eigenvalues of the matrix
and (as we will show) the symmetric group
characters define a basis of inhomogeneous
degree for the space of symmetric functions.

To arrive at an answer to this question and its role in representation
theory we begin in the
context of Schur-Weyl duality.  Consider an $n$-dimensional
vector space $V = \{ v_1, v_2, \ldots, v_n\}$ for which
$Gl_n$ has an action on $V$ and 
also acts diagonally on the tensor
space $T^k(V)$.  The symmetric group $S_k$ acts by permutingthe positions
of the tensor space and the symmetric group algebra is the centralizer
algebra of the $Gl_n$ action.  Schur-Weyl duality states that $T^k(V)$
decomposes as a direct sum of irreducible $S_k$ and $Gl_n$ bimodules
\begin{equation}
T^k(V) \simeq \bigoplus_{\lambda \vdash k} {\bf S}^\lambda \otimes M^\lambda
\end{equation}
where the ${\bf S}^\lambda$ is an irreducible $S_k$ module
and $M^\lambda$ is an irreducible $Gl_n$ module.  The Schur functions
arise in two different contexts.  The first as the
Frobenius image of the character of ${\bf S}^\lambda$,
\begin{equation}
s_\lambda = \frac{1}{k!} \sum_{\sigma \in S_k} char_{{\bf S}^\lambda}(\sigma) p_{type(\sigma)},
\end{equation}
and the second as characters
\begin{equation}
s_\lambda(x_1, x_2, \ldots, x_n) = char_{M^\lambda}(A)
\end{equation}
where $A$ is a matrix whose eigenvalues are $x_1, x_2, \ldots, x_n$.

It is this second context in which we propose to consider 
the irreducible characters
of the symmetric group $S_n$ which sits inside of 
$Gl_n$ as permutation matrices.
We wish to consider the question that we posed
in the first paragraph where we change the group $Gl_n$ with
another group.  In the case that the group $Gl_n$ is
replaced with $SO_{2n}$, $SO_{2n+1}$ or $Sp_{2n}$, formulas
for the characters were given by Weyl \cite{Weyl} and
explicit expressions of 
the characters as symmetric functions were
worked out by Koike and Terada \cite{KT}.
The centralizer algebras of the group action in these cases
are the Brauer algebras.

It is the case when group $Gl_n$
is replaced by the symmetric group $S_n$ of permutation matrices
that we are considering here.  The centralizer algebra of this action was
determined by Jones and Martin \cite{Jo,Ma1,Ma2,Ma3,Ma4} to be the partition algebra.
Halverson and Ram \cite{Halverson, HalversonRam} made
the characters of the partition algebra explicit and we see
these character values arise in the coefficients of the change of
bases.  We have not however seen the symmetric group characters
treated in a manner similar to the characters of the classical groups.

To this end, we propose the introduction of the
{\it irreducible symmetric group characters} as a basis
(or, for short, irreducible character basis), 
$\st_{\lambda}$, of the symmetric
functions such that when this function is evaluated
at the eigenvalues of a permutation matrix for the permutation $\sigma$,
the value is equal to the irreducible character
$\chi^{(|\sigma|-|\la|,\lambda)}(\sigma)$.  That the character only
depends on the size of the permutation and cycle structure is a well known
fact of characters of the symmetric group, but it will also be a consequence
of the way that we define the symmetric function $\st_\lambda$.

To understand the basic properties of this symmetric function
and the change of basis coefficients, a natural second basis is suggested
by the algebra. We call this second basis the {\it induced
trivial character basis}, $\hht_\la$, and these functions are equal to
the characters of the induced trivial representation from
from $S_{\la_1} \times S_{\la_2} \times \cdots \times S_{\la_{\ell(\la)}}$
to $S_n$ where $n$ is some positive integer.  Like the
irreducible character basis, the induced trivial character
$\hht_\lambda$ is also independent of
the value of $n$.

Because $\st_\la$ and $\hht_\la$ are characters, their products
are the characters of the tensors of the representations with
the diagonal action.
The structure coefficients of the irreducible character
basis will be the {\it reduced Kronecker coefficients} \cite{BOR}.
That is, for an $n$ sufficiently large and $\la, \mu$ partitions
of $n$, for any partition $\nu$ let ${\overline \nu} = (\nu_2, \nu_3,
\ldots, \nu_{\ell(\nu)})$, then
\begin{equation}
\st_{\overline \la} \st_{\overline \mu} = \sum_{\nu \vdash n}
k_{\la\mu\nu} \st_{\overline \nu}
\hskip .3in
\hht_{\overline \la} \hht_{\overline \mu} = \sum_{\nu \vdash n}
{\mathcal K}_{\la\mu\nu} \hht_{\overline \nu}
\end{equation}
where for partitions $\alpha$, $\beta$, $\tau$ of size $n$, and
\begin{equation}
k_{\alpha\beta\tau} =
\sum_{\ga \vdash n}
\chi^{\alpha}(\ga)
\chi^{\beta}(\ga)
\chi^{\tau}(\ga)
\hbox{ and }
{\mathcal K}_{\alpha\beta\tau} =
\sum_{\ga \vdash n}
\left< h_\alpha, p_\ga \right>
\left< h_\beta, p_\ga \right>
\left< h_\tau, p_\ga \right>~.
\end{equation}
The $k_{\alpha\beta\tau}$
and ${\mathcal K}_{\alpha\beta\tau}$
are the coefficients that appear in the internal product
of the Schur and complete symmetric function bases.

In addition, due to the Schur-Weyl duality with the partition algebra,
the change of basis coefficients between the power sum basis
$p_\mu$ and the elements $\st_\lambda$ are the irreducible characters
of the partition algebra and our elements must satisfy the
Murgnahan-Nakayma rule for those characters (see \cite{Halverson}). In addition, the product $\st_{\overline \la} \st_{\overline \mu}$ corresponds to the restriction of characters of the partition algebra \cite{BDO}.

A recent paper by Church and Farb \cite{ChurchFarb} introduces
a notion of representation stability.  They describe a number of
families of representations whose decomposition into irreducible
representations is independent of $n$, if $n$ is sufficiently large.
The bases that we introduce in this paper are likely to be
a useful tool in finding expressions for their character.

There are a number of other results that we are not adequately
able to summarize in this abstract.  Some of these include
product rules which come from known cases of Kronecker coefficients,
change of basis coefficients, coproducts and connections with
other questions in representation theory.

We believe that the relationship of the representation theory,
combinatorics and symmetric functions already merits study of
the induced character and irreducible character bases, however
we are aware that these bases will only be seen as an
important innovation in the theory of symmetric functions
if they find proper applications.

Since inner and outer plethysm and the inner (Kronecker) product 
of Schur functions are
all classes of symmetric function expressions which do not currently
have a simple combinatorial formula, and yet they all three have
the a stability property that is related to removing the first
row, we expect that this basis will prove a useful tool
for progress in explaining these operations.
\end{section}

\begin{section}{Notation}
The objects that arise in
the following constructions are familar
building blocks of combinatorics: set, multi-set, partition,
set partition, multi-set partition,
composition, weak composition,
set composition, multi-set composition,
tableau, words, etc.  
We will use this section to establish
notation conventions and use of language
that we will need in this paper.

For non-negative integers $n$ and $\ell$, a partition
of size $n$ and length $\ell$ is a sequence
$\lambda = (\lambda_1, \lambda_2, \ldots, \lambda_{\ell})$
such that $\lambda_i \geq \lambda_{i+1}$ for $1 \leq i < \ell$
and $\lambda_1 + \lambda_2 + \cdots + \lambda_\ell = n$.
The size of the partition is denoted $|\la| = n$ and
the length of the partition is denoted $\ell(\la)=\ell$.
We will often use the shorthand notation $\la \vdash n$
to indicate that $\lambda$ is a partition of $n$.
The symbols $\lambda$ and $\mu$ will be reserved exclusively
for partitions.  Let $m_i(\la)$ represent the number of
times that $i$ appears in the partition $\la$.  It will be
convenient to sometimes represent our partitions in
exponential notation where $m_i = m_i(\la)$ and
$\la = (1^{m_1}2^{m_2}\cdots k^{m_k})$. With this notation
the number of permutations with cycle structure $\lambda\vdash n$ is
$\frac{n!}{z_\la}$ where
\begin{equation}
z_\la = \prod_{i=1}^{\la_1} m_i(\la)! i^{m_i(\la)}~.
\end{equation}
The most common operation we will use will be adding a
part of size $n-|\la|$ to the beginning of a partition.  This
will be denoted $(n-|\la|, \la)$.  If $n < |\la|+\la_1$, this
sequence will no longer be an integer partition.

The cells of a partition $\la$ are the set of points
$\{ (i,j) : 1 \leq i \leq \la_j, 1 \leq j \leq \ell(\la) \}$.
We will represent these cells as stacks of boxes in
the first quadrant (following `French notation' for a
partition).  A tableau is a mapping from the set of
cells to a set of labels and a tableau will be represented
by filling the boxes of the diagram for a partition with
the labels.
In our case, we will
encounter tableaux where only a subset of the cells are
mapped to a label.

Multi-sets will also be represented by
exponential notation so that
$\{ 1^{a_1},2^{a_2},\ldots, \ell^{a_\ell}\}$
represents the multi-set where $i$ occurs
$a_i$ times.

A set partition of a set $S$ is a set of subsets
$\{ S_1, S_2, \ldots, S_\ell\}$ with
$S_i \cap S_j = \emptyset$ for $1 \leq i< j \leq \ell$
and
$S_1 \cup S_2 \cup \cdots \cup S_\ell = S$.
A multi-set partition $\pi = \{S_1, S_2, \ldots, S_\ell\}$ 
of a multi-set $S$ is a similar construction to a
set partition, but now $S_i$ may be
a multi-set, and it is possible that two sets
$S_i$ and $S_j$ have non-empty intersection (and may
even be equal).  Let $\ell(\pi) := \ell$ and
$\mt(\pi)$ represent the partition of $\ell(\pi)$
consisting of
the the multipicities of the multi-sets which
occur in $\pi$ (e.g. $\mt(\{\{1,1,2\},\{1,1,2\},\{1,3\}\})
=(2,1)$ because $\{1,1,2\}$ occurs 2 times and 
$\{1,3\}$ occurs 1 time).  We will use the notation
$\pi \mvdash S$ to indicate that $\pi$ is a multi-set
partition of the multi-set $S$.

%

For non-negative integers $n$ and $\ell$, a composition 
is an ordered sequence 
$\alpha = (\alpha_1, \al_2, \ldots, \al_\ell)$
of integers $\al_i > 0$ such that
$\al_1+\al_2 + \cdots +\al_\ell=n$.  A weak composition
is such a sequence with the condition that $\al_i \geq 0$.
To indicate that $\alpha$ is a composition of $n$ we will
use the notation  
$\alpha \models n$ and to indicate that
$\alpha$ is a weak composition of $n$ we will use the notation
$\alpha \models_w n$.  For both
compositions and weak compositions, $\ell(\alpha) := \ell$.

The ring of symmetric functions
(for some modern references on this subject see
for example \cite{Mac, Sagan, Stanley, Lascoux})
will be denoted
$Sym = \QQ[p_1, p_2, p_3, \ldots]$ and has
the fundamental bases (each indexed by the set of
partitions $\lambda$)
{\it power sum} $\{p_\la\}_{\la}$,
{\it homogeneous/complete} $\{h_\la\}_{\la}$,
{\it elementary} $\{e_\la\}_{\la}$,
{\it Schur} $\{s_\la\}_{\la}$
and {\it monomial} $\{m_\la\}_{\la}$.
We will also refer to the irreducible character
of the symmetric group indexed by the partition $\la$
and evaluated at a permutation
of cycle structure $\mu$ as the coefficient
$\left< s_\la, p_\mu \right> = \chi^\la(\mu)$.

For $k > 0$, define
$$\Xi_k := 1, e^{2\pi i/k}, e^{4\pi i/k},
\ldots, e^{2(k-1)\pi i/k}$$
as a symbol representing the eigenvalues of a permutation
matrix of a $k$-cycle.  Then for any partition
$\mu$, let
$$\Xi_\mu := \Xi_{\mu_1}, \Xi_{\mu_2}, 
\ldots, \Xi_{\mu_{\ell(\mu)}}$$
be the multi-set of eigenvalues of a permutation matrix with
cycle structure $\mu$.  We will evaluate symmetric
functions at these eigenvalues.  The notation
$f[\Xi_\mu]$ represents taking the element
$f \in Sym$ and replacing $p_k$ in $f$ with
$x_1^k + x_2^k + \cdots + x_{|\mu|}^k$ and then
replacing the variables $x_i$ with the values in $\Xi_\mu$.
\end{section}

\begin{section}{Symmetric group character bases 
of the symmetric functions}

To begin, we establish the following result
connecting the evaluation of the homogeneous symmetric functions
at roots of unity $\Xi_\mu$ with the values of the trivial character
induced from $S_{\nu_1} \times S_{\nu_2} \times \cdots
\times S_{\nu_{\ell(\nu)}}$ to $S_{|\mu|}$.  This is given
by the scalar product $\left< h_{\nu} h_{|\mu|-|\nu|}, p_\mu \right>$.

\begin{theorem}
For all partitions $\lambda$,
\begin{equation}\label{eq:hlainhtbasis}
h_\lambda[\Xi_\mu] = 
\sum_{\pi \mvdash \{1^{\la_1},2^{\la_2},\ldots,\ell^{\la_\ell}\}}
\left< h_{\mt(\pi)} h_{|\mu|-|\mt(\pi)|}, p_\mu \right>.
\end{equation}
\end{theorem}

We may also extend this further by computing a similar
expression for the elementary basis.  While we don't include this
result here as the result is auxiliary to what we are presenting,
the expansion
involves combinatorial objects which are set partitions of a
multi-set (as opposed to multi-set partitions of a multi-set)
since in each set of the set partition, repeated elements are
not allowed.

We are now ready to define symmetric functions whose
evaluations at roots of unity (the eigenvalues of a
permutation matrix) are the values of characters.
We base our definition on equation \eqref{eq:hlainhtbasis}.

\begin{definition}\label{def:hht}
Let $\hht_\mu$ be the family of symmetric functions which
satisfies:

\begin{equation}\label{eq:definingrel}
h_\lambda = 
\sum_{\pi \mvdash \{1^{\la_1},2^{\la_2},\ldots,\ell^{\la_\ell}\}}
\hht_{\mt(\pi)}.
\end{equation}
We call this basis the {\it induced trivial character} basis of
the symmetric functions.
\end{definition}
This is a recursive definition for calculating this basis
directly since there is precisely one multi-set partition
of $\{1^{\la_1},2^{\la_2},\ldots,\ell^{\la_\ell}\}$
such that $\mt(\pi)$ is of size $|\la|$, hence

\begin{equation}
\hht_\la = h_\la - \sum_{\substack{\pi \mvdash \{1^{\la_1},2^{\la_2},\ldots,\ell^{\la_\ell}\}\\\mt(\pi)\neq\la}}
\hht_{\mt(\pi)}.
\end{equation}

Now by equation \eqref{eq:hlainhtbasis} and an induction
argument, we can conclude that for all partitions $\mu$,
\begin{equation}
\hht_{\la}[\Xi_\mu] = \left< h_{|\mu|-|\la|}h_{\la}, p_\mu\right>
\end{equation}
and this is the value of the character of the
trivial module which is induced from 
$S_{|\mu|-|\la|} \times S_{\la_1} \times S_{\la_2} \times
\cdots \times S_{\la_{\ell(\la)}}$ to the full symmetric group
$S_{|\mu|}$.

\begin{example}
We list all of the multi-set partitions of $\{1,1,1,2\}$.  They
are
$$\{\{1\},\{1\},\{1\},\{2\}\}, \{\{1\},\{1\},\{1,2\}\},
\{\{1,1\},\{1\},\{2\}\},\{\{1,1,1\},\{2\}\},$$
$$\{\{1,1,2\},\{1\}\},\{\{1,1\},\{1,2\}\},\{\{1,1,1,2\}\}$$
By Definition \ref{def:hht},
$$h_{31} = \hht_{31} + \hht_{21} + \hht_{111} + 3\hht_{11} + \hht_1.$$
\end{example}

Next we define the symmetric functions $\st_\la$ by the
change of basis with $\hht_\la$ basis and the formulas
\begin{equation}\label{eq:httost}
\hht_\mu = \sum_{|\la|\leq|\mu|} K_{(n-|\la|,\la)(n-|\mu|,\mu)} \st_\la
\end{equation}
and
\begin{equation}\label{eq:sttoht}
\st_\la = \sum_{|\mu|\leq|\la|} 
K_{(n-|\la|,\la)(n-|\mu|,\mu)}^{-1} \hht_\mu
\end{equation}
where $n$ is any positive integer greater than $max(|\la|+\la_1,|\mu|+\mu_1)$ and
$K_{\la\mu}$ are the Kostka coefficients (the
change of basis coefficients between the complete
symmetric functions and the Schur basis, or the
number of column strict tableaux of shape $\la$
and content $\mu$).  It should be clear from a
number of perspectives that this
definition is independent of the value of $n$
as long as $n$ is sufficiently large.
In particular, if $n-|\la|\geq \la_1$
then the Kostka coefficients $K_{(n-|\la|,\la)(n-|\mu|,\mu)}$
do not change by increasing the value of $n$ since there
is an isomorphism between
the tableau with $n-|\mu|$ labels $1$ in the first row
and those that have $n-|\mu|+1$ labels $1$ in the first row
where the length of the first row also differs by $1$.

If $n$ is smaller than $|\la|+\la_1$, then the
change of basis coefficients are the same as those
between the complete symmetric functions and
a Schur function indexed by a composition
$\alpha = (|\mu|-|\la|,\la)$, namely
the expression representing the Jacobi-Trudi
matrix
\begin{equation}
\det\left[ h_{\alpha_i + i - j} \right]_{1 \leq i,j \leq \ell(\la)+1}.
\end{equation}

We find then that $\st_\la$ are the (unique) symmetric
functions of inhomogeneous degree $|\la|$, such that
\begin{align*}
\st_\la[\Xi_\ga]
&= \sum_{|\mu|\leq|\la|} 
K_{(n-|\la|,\la)(n-|\mu|,\mu)}^{-1}
\hht_\mu[\Xi_\ga]\\
&=
\sum_{|\mu|\leq|\la|} 
K_{(n-|\la|,\la)(n-|\mu|,\mu)}^{-1}
\left< h_{(n-|\mu|,\mu)}, p_\ga \right>\\
&= \left<s_{(n-|\la|,\la)}, p_\ga\right>
= \chi^{(n-|\la|,\la)}(\ga).
\end{align*}

We conclude that the basis $\st_\la$ are the characters of the
symmetric group when the symmetric group is realized as
permutation matrices acting on the representation.

The combinatorial interpretation of the coefficients for the
$\st$-expansion of a complete symmetric function comes from
combining the notion of multi-set partition of a multi-set
and column strict tableau.
\begin{prop} \label{prop:stexpofh}
For a partition, choose an $n\geq\mu_1+|\mu|$, then
\begin{equation}
h_\mu = \sum_{T}
\st_{\overline {shape(T)}}
\end{equation}
where ${\overline \la} = (\la_2, \la_3,\ldots, \la_{\ell(\la)})$ and
the sum is over all column strict tableau with $n - |\mu|$
blank cells in the first row and the rest of the cells filled
with multi-sets of labels such that the total content of the
tableau is $\{1^{\mu_1}2^{\mu_2}\cdots\ell^{\mu_\ell}\}$.
\end{prop}

\begin{example}
Consider the following 20 column strict tableaux whose
entries are multi-sets with total content of the tableau
$\{1^2,2\}$ and with $5$ blank cells.

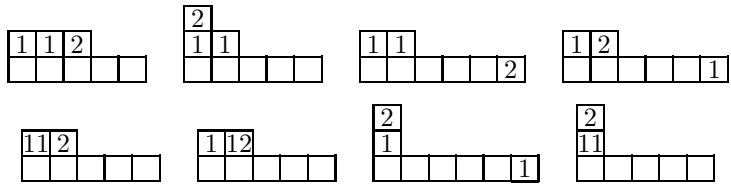
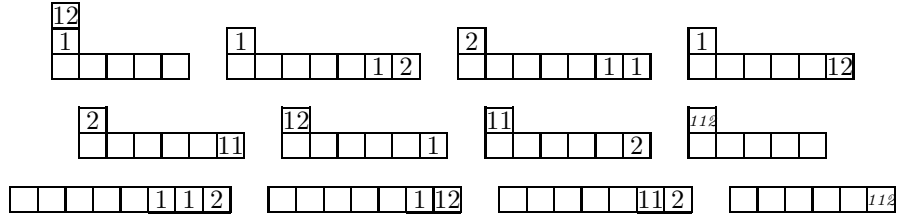
\begin{equation*}
\young{1&1&2\cr&&&&\cr}\hskip .2in
\young{2\cr1&1\cr&&&&\cr}\hskip .2in
\young{1&1\cr&&&&&2\cr}\hskip .2in
\young{1&2\cr&&&&&1\cr}
\end{equation*}
\begin{equation*}
\young{11&2\cr&&&&\cr}\hskip .2in
\young{1&12\cr&&&&\cr}\hskip .2in
\young{2\cr1\cr&&&&&1\cr}\hskip .2in
\young{2\cr11\cr&&&&\cr}
\end{equation*}
\begin{equation*}
\young{12\cr1\cr&&&&\cr}\hskip .2in
\young{1\cr&&&&&1&2\cr}\hskip .2in
\young{2\cr&&&&&1&1\cr}\hskip .2in
\young{1\cr&&&&&12\cr}
\end{equation*}
\begin{equation*}
\young{2\cr&&&&&11\cr}\hskip .2in
\young{12\cr&&&&&1\cr}\hskip .2in
\young{11\cr&&&&&2\cr}\hskip .2in
\young{\hbox{\tiny 112}\cr&&&&\cr}
\end{equation*}
\begin{equation*}
\young{&&&&&1&1&2\cr}\hskip .2in
\young{&&&&&1&12\cr}\hskip .2in
\young{&&&&&11&2\cr}\hskip .2in
\young{&&&&&\hbox{\tiny 112}\cr}
\end{equation*}
Proposition \ref{prop:stexpofh} states then that
\begin{equation}
h_{21} = 4\st_{()} + 7\st_1 + 3\st_{11} + 4\st_2 + \st_{21} + \st_3~.
\end{equation}
\end{example}


\begin{example}\label{ex:settableaux}
Say that we want to compute the decomposition of $V^{\otimes 4}$
where $V = {\mathcal L}\{ x_1, x_2, x_3, \ldots ,x_n \}$ as an $S_n$ module
with the diagonal action.  The module $V$ has character equal to
$\hht_1 = h_1$.  Therefore to compute the decomposition of this
character into $S_n$ irreducibles we are looking for the expansion
of $h_{1^4}$ into the $\st$-basis.

Using Sage we compute that it is
\begin{align}
h_{1^4} = &15\st_{()} + 37\st_{1} + 31\st_{11} + 10\st_{111} + \st_{1111} + 31\st_{2}\\
&+ 20\st_{21} + 3\st_{211} + 2\st_{22} + 10\st_{3} + 3\st_{31} + \st_{4}
\end{align}

If $n\geq6$ then the multiplicity of the irreducible $(n-3,3)$ will be $10$.
The combinatorial interpretation of this value is the number of column strict
tableaux with entries that
are multi-sets (or in this case sets) of $\{1,2,3,4\}$ of shape $(4,3)$ or $(5,3)$
and $4$ blank entries in the first row. Those tableaux are
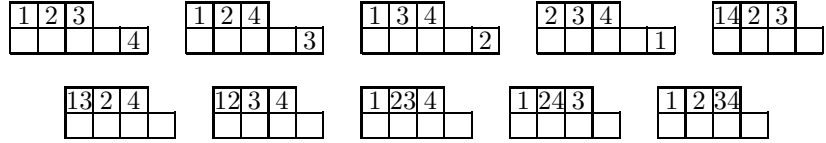
\begin{equation*}
\young{1&2&3\cr&&&&4\cr}\hskip .2in
\young{1&2&4\cr&&&&3\cr}\hskip .2in
\young{1&3&4\cr&&&&2\cr}\hskip .2in
\young{2&3&4\cr&&&&1\cr}\hskip .2in
\young{14&2&3\cr&&&\cr}
\end{equation*}
\begin{equation*}
\young{13&2&4\cr&&&\cr}\hskip .2in
\young{12&3&4\cr&&&\cr}\hskip .2in
\young{1&23&4\cr&&&\cr}\hskip .2in
\young{1&24&3\cr&&&\cr}\hskip .2in
\young{1&2&34\cr&&&\cr}
\end{equation*}
What is interesting about this example is that the usual combinatorial
interpretation for the repeated Kronecker product $(\chi^{(n-1,1)})^{\ast k}$ is
stated in other places in the literature in terms of oscillating tableaux
\cite{CG, Sund}
so this more general formula is giving a slightly different way of combinatorially
describing the multiplicities in terms of set valued tableaux.
\end{example}
\end{section}

\begin{section}{Character polynomials and the irreducible
character basis}

Character polynomials were first used
by Murnaghan \cite{Murg}. Much later, Specht
\cite{Sp} gave determinantal formulas and expressions in
terms of binomial coefficients for these polynomials.
They are treated as an example
in Macdonald's book \cite[ex. I.7.13 and I.7.14]{Mac}.
More recently, Garsia and Goupil \cite{GG} 
gave an umbral formula for computing them.  We will
show in this section that character polynomials
are a transformation of character symmetric functions
and this will allow us to give an expression
for character symmetric functions in the power sum basis.

To begin, we note that
$p_r[\Xi_k] = k$ if $k$ divides $r$ and it is equal to $0$
otherwise.
In general, we can express any partition
$\mu$ in exponential notation 
$\mu = (1^{m_1} 2^{m_2}\cdots r^{m_r})$
where $m_i$ are the number of parts of size $i$ in $\mu$.
Therefore
\begin{equation}
p_k[\Xi_\mu] = \sum_{d|k} d m_d
\end{equation}
Hence any symmetric function $f$ evaluated at
some set of roots of unity is equal to 
a polynomial in variables $m_1, m_2, \ldots, m_n$
where
\begin{equation} \label{eq:sftopoly}
f[\Xi_\mu] = f \Big|_{p_k \rightarrow \sum_{d|k} d m_d}
= q( m_1, m_2, \ldots, m_n )~.
\end{equation}
Moreover, if we know this polynomial $q(m_1, m_2, \ldots, m_d)$
we can use M\"obius inversion to recover
the symmetric function since if $p_k = \sum_{d|k} d m_d$,
then $k m_k = \sum_{d|k} \mu(k/d) p_d$ where
\begin{equation}
\mu(r) = \begin{cases}(-1)^d&\hbox{if $r$ is a
product of $d$ distinct primes}\\
0&\hbox{if $r$ is not square free}
\end{cases}.
\end{equation}
Therefore, we also have
\begin{equation}\label{eq:recoversf}
q\left( p_1, \frac{p_2-p_1}{2},\frac{p_3-p_1}{3}, 
\ldots, \frac{1}{n}\sum_{d|n} \mu(n/d) p_d \right) = f~.
\end{equation}

Following the notation of \cite{GG},
a character polynomial is a multivariate
polynomial $q_\lambda(x_1, x_2, x_3, \ldots)$ in the variables
$x_i$ such that for specific integer values $x_i = m_i \in {\mathbb Z}$,
\begin{equation}
q_\lambda(m_1, m_2, m_3, \ldots)
= \chi^{(n-|\lambda|,\lambda)}(1^{m_1}2^{m_2}3^{m_3}\cdots)
\end{equation}
where $n = \sum_{i \geq 1} i m_i$.  As a consequence
of this relationship we have the following
relationship between the character polynomials
$q_\lambda( x_1, x_2, x_3, \ldots)$
and character basis $\st_\la$.

\begin{prop} For a partition $\la$,
$$q_\lambda(x_1, x_2, x_3, \ldots) = \st_\la \Big|_{p_k \rightarrow \sum_{d|k} d x_d}$$
and
$$\st_\la = q_\lambda(x_1, x_2, x_3, \ldots) \Big|_{x_k \rightarrow \frac{1}{k}\sum_{d|k} \mu(k/d) p_d }$$
\end{prop}

As an important intermediate result, we have
a power sum expansion of the irreducible character basis.
Consider the following example.

\begin{example}  In \cite{GG}, the formula that they use to
compute the character polynomial is most simply stated algorithmically
at the top of page 3.  If we make an additional substitution
in the last step of their formula with $x_k$ by
$\frac{1}{k} \sum_{d|k} \mu(k/d) p_d$, then
in order to compute $\st_\la$ we compute the four following steps,
plus a fifth step to relate the character polynomial to
irreducible character symmetric function.
\begin{enumerate}
\item Expand the Schur function $s_\lambda$ in the power sums basis
$s_\la = \sum_{\alpha} \frac{\chi^\lambda(\ga)}{z_\ga} p_\ga$.
\item Replace each power sum $p_i$ by $i x_i -1$.
\item Expand each product $\prod_i (ix_i-1)^{a_i}$ 
as a sum $\sum_g c_g \prod_i x_i^{g_i}$.
\item Replace each $x_k^{g_k}$ by $(x_k)_{g_k} =
x_k(x_k-1)\cdots(x_k-g_k+1)$.
\item Replace each $x_k$ by $\frac{1}{k} \sum_{d|k} \mu(k/d) p_d$.
\end{enumerate}
For instance, to compute $\st_3$ we follow the steps to obtain:
\begin{enumerate}
\item $s_3 = \frac{1}{6}(p_{1^3} + 2 p_{21} + 3 p_3)$
\item $\frac{1}{6}(p_{1^3} + 2 p_{21} + 3 p_3) \rightarrow
\frac{1}{6}((x_1-1)^3+3(2x_2-1)(x_1-1)+2(3x_3-1))$
\item $\frac{1}{6}((x_1-1)^3+3(2x_2-1)(x_1-1)+2(3x_3-1))=
\frac{1}{6}x_1^3-\frac{1}{2}x_1^2+x_1x_2-x_2+x_3$
\item $q_{3} = \frac{1}{6}(x_1)_3 - \frac{1}{2}(x_1)_2
+ x_1x_2 - x_2+ x_3$
\item $\st_{3} = \frac{1}{6}(p_1)_3 - \frac{1}{2}(p_1)_2
+ p_1\frac{p_2-p_1}{2} - \frac{p_2-p_1}{2}+ \frac{p_3-p_1}{3}$
\end{enumerate}
\end{example}

As a consequence we state a slightly more explicit
expression for the character basis that we derive by
following the algorithm of \cite{GG}. 

\begin{prop} \label{prop:pexpansionofst} For $\lambda \vdash n$,
\begin{equation}\label{eq:pexpansionofst}
\st_\lambda = \sum_{\gamma \vdash n} 
\chi^{\lambda}(\gamma) \frac{\bfp_{\ga}}{z_\gamma}
\end{equation}
where
\begin{equation}\label{eq:weirdp}
\bfp_{i^r} = \sum_{k=0}^r (-1)^{r-k} i^k \binom{r}{k}
\left(\frac{1}{i} \sum_{d|i} \mu(i/d) p_d \right)_k\hbox{ and }\bfp_{\gamma}
:= \prod_{i \geq 1} {\mathbf p}_{i^{m_i(\gamma)}}~.
\end{equation}
\end{prop}
\end{section}


\begin{section}{A product rule for the $\hht$-basis}
The Kronecker product of two complete symmetric functions
expanded in the complete basis is listed as Exercise 23 (e) in section I.7
of \cite{Mac}.
It states that
\begin{equation}\label{eq:usualhproduct}
h_\la \ast h_\mu = \sum_{M} \prod_{i=1}^{\ell(\la)} \prod_{j=1}^{\ell(\mu)} h_{M_{ij}}
\end{equation}
summed over all matrices $M$ of non-negative integers with $\ell(\la)$
rows, $\ell(\mu)$ columns and row sums $\la_i$ and column sums $\mu_j$.

We propose a different (but equivalent) 
combinatorial interpretation for this product based on the $\hht$-basis.
For two multi-set partitions $\pi = \{ S_1, S_2, \ldots, S_{\ell(\pi)}\}$
and $\theta = \{T_1, T_2, \ldots, T_{\ell(\theta)}\}$ with $S_i \cap T_j = \emptyset$.
Define a join operation (which we will denote as $\pi\odot\tau$)
as the set of distinct multi-set partitions of the form
$$\{S_{i_1}, S_{i_2}, \ldots, S_{i_{\ell(\pi)-r}}, 
T_{j_1}, T_{j_2}, \ldots, T_{\ell(\theta)-r},
S_{u_1} \cup T_{v_1}, S_{u_2} \cup T_{v_2}, \ldots, 
S_{u_r} \cup T_{v_r}\}$$
with $\{i_1, i_2, \ldots, i_{\ell(\pi)-r}, u_1, u_2, \ldots, u_r\}
= \{1,2,\ldots,\ell(\pi)\}$
and $\{j_1, j_2, \ldots, j_{\ell(\theta)-r}, v_1, v_2, \ldots, v_r\}
= \{1,$ $2, \ldots, \ell(\theta)\}$.

\begin{prop}\label{prodht} For disjoint multi-set partitions $\pi$ and $\theta$,
$$\hht_{\mt(\pi)} \hht_{\mt(\theta)} = \sum_{\tau \in \pi\odot\tau} \hht_{\mt(\tau)}.$$
\end{prop}

\begin{example}
let $\pi = \{\{1\},\{1\},\{2\}\}$ and $\theta = \{\{3\},\{3\},\{4\}\}$

Below we list the multi-set partitions in the product
along with the corresponding partition $\mt(\tau)$.
$$\{\{1\},\{1\},\{2\},\{3\},\{3\},\{4\}\} \rightarrow (2,2,1,1)
\hskip .3in\{\{1,3\},\{1\},\{2\},\{3\},\{4\}\} \rightarrow (1,1,1,1,1)$$
$$\{\{1\},\{1\},\{2,3\},\{3\},\{4\}\} \rightarrow (2,1,1,1)
\hskip .3in\{\{1\},\{1,4\},\{2\},\{3\},\{3\}\} \rightarrow (2,1,1,1)$$
$$\{\{1\},\{1\},\{2,4\},\{3\},\{3\}\} \rightarrow (2,2,1)
\hskip .3in\{\{1,3\},\{1,3\},\{2\},\{4\}\} \rightarrow (2,1,1) $$
$$\{\{1,3\},\{1,4\},\{2\},\{3\}\} \rightarrow (1,1,1,1)
\hskip .3in\{\{1\},\{1,3\},\{2,3\},\{4\}\} \rightarrow (1,1,1,1)$$
$$\{\{1\},\{1,4\},\{2,3\},\{3\}\} \rightarrow (1,1,1,1)
\hskip .3in\{\{1\},\{1,3\},\{2,4\},\{3\}\} \rightarrow (1,1,1,1)$$
$$\{\{1,3\},\{1,3\},\{2,4\}\} \rightarrow (2,1)
\hskip .3in\{\{1,3\},\{1,4\},\{2,3\}\} \rightarrow (1,1,1)$$
As a consequence of Proposition \ref{prodht} we conclude
\begin{equation}
\hht_{21} \hht_{21} =
\hht_{111} + 4\hht_{1111} + \hht_{11111}
+ \hht_{21} + \hht_{211} + 2\hht_{2111} + \hht_{221}
+ \hht_{2211}
\end{equation}
or in terms of Kronecker products
\begin{equation}
h_{521} \ast h_{521} =
h_{5111} + 4h_{41111} + h_{311111}
+ h_{521} + h_{4211} + 2h_{32111} + h_{3221}
+ h_{22211}~.
\end{equation}
\end{example}
\end{section}

\begin{section}{Appendix: Using Sage to compute
the character bases of symmetric functions}

Sage is an open source symbolic calculation program based on
the computer language Python.  A large community of mathematicians
participate in its support and add to its functionality
\cite{sage, sage-co}.  In particular, the built-in library
for symmetric functions includes a large extensible set of functions
which makes it possible to do calculations within the ring
following closely the mathematical notation that we use in this paper.
The language itself has a learning curve, but the contributions
made by the community towards the functionality
make that a barrier worth overcoming.

In version 6.10 or later of Sage (released Jan 2016)
these bases will be available
as methods in the ring of symmetric functions.

We demonstrate examples of some of the definitions and results in this paper
by sage calculations.
\begin{verbatim}
sage: Sym = SymmetricFunctions(QQ)
sage: Sym
Symmetric Functions over Rational Field
sage: st = Sym.irreducible_symmetric_group_character()
sage: st
Symmetric Functions over Rational Field in the irreducible character basis
\end{verbatim}
We can compare the structure coefficients of the $\st$-basis with
the Kronecker product of Schur functions whose first part is
sufficiently large.
\begin{verbatim}
sage: st[2]*st[2]
st[] + st[1] + st[1, 1] + st[1, 1, 1] + 2*st[2] + 2*st[2, 1] + st[2, 2]
+ st[3] + st[3, 1] + st[4]
sage: s = Sym.Schur(); s
Symmetric Functions over Rational Field in the Schur basis
sage: s[6,2].kronecker_product(s[6,2])
s[4, 2, 2] + s[4, 3, 1] + s[4, 4] + s[5, 1, 1, 1] + 2*s[5, 2, 1] + s[5,
3] + s[6, 1, 1] + 2*s[6, 2] + s[7, 1] + s[8]
\end{verbatim}

We can express one basis in terms of another.
\begin{verbatim}
sage: h = Sym.complete(); h
sage: Symmetric Functions over Rational Field in the homogeneous basis
sage: ht = Sym.induced_trivial_character(); ht
Symmetric Functions over Rational Field in the induced trivial character basis
sage: ht(h[2,2]) # express h_22 in the ht-basis
ht[1] + 3*ht[1, 1] + ht[1, 1, 1] + ht[2] + 2*ht[2, 1] + ht[2, 2]
sage: st(h[2,2])
9*st[] + 17*st[1] + 9*st[1, 1] + st[1, 1, 1] + 14*st[2] + 6*st[2, 1] +
st[2, 2] + 5*st[3] + st[3, 1] + st[4]
sage: st(ht[2,2])
st[] + 2*st[1] + st[1, 1] + 3*st[2] + 2*st[2, 1] + st[2, 2] + 2*st[3] +
st[3, 1] + st[4]
sage: s(h[4,2,2])
s[4, 2, 2] + s[4, 3, 1] + s[4, 4] + 2*s[5, 2, 1] + 2*s[5, 3] + s[6, 1,
1] + 3*s[6, 2] + 2*s[7, 1] + s[8]
\end{verbatim}
This ticket also introduces two operations that we use in this paper.
The first is the evaluation of a symmetric function at the eigenvalues of
a permutation matrix where the permutation has cycle structure $\mu$.
We represented this operation in the paper as $f[\Xi_\mu]$.
In Sage, elements of the symmetric functions have the method {\tt eval\_at\_permutation\_roots}
which represents this operation.
\begin{verbatim}
sage: ht[3,1].eval_at_permutation_roots([3,3,2,2,1])
2
sage: st[3,1].eval_at_permutation_roots([3,3,2,2,1])
-1
\end{verbatim}
In Sage, these can be compared to the following coefficients
in the power sum basis.
\begin{verbatim}
sage: p = Sym.powersum(); p
Symmetric Functions over Rational Field in the powersum basis
sage: h[7,3,1].scalar(p[3,3,2,2,1])
2
sage: s[7,3,1].scalar(p[3,3,2,2,1])
-1
\end{verbatim}

The other operation that we define is one that
interprets a symmetric function
as a symmetric group character and then
maps that character
to the Frobenius image (or characteristic map)
of the character.
That is, for a symmetric function $f$,
the function computes
\begin{equation}
\phi_n(f) = \sum_{\mu \vdash n} f[\Xi_\mu] \frac{p_\mu}{z_\mu}~.
\end{equation}
The elements of the symmetric functions in Sage
will have the method {\tt character\_to\_frobenius\_image} which
represents the map $\phi_n$.

We have defined the $\st$ and $\hht$ bases so that they have the property
$\phi_n(\st_\lambda) = s_{(n-|\la|,\la)}$ and $\phi_n(\hht_\lambda) = h_{(n-|\la|,\la)}$ respectively
if $n \geq |\la|+\la_1$.  If $n<|\la|+\la_1$, then the
corresponding symmetric function will still be equivalent to 
a Schur function or complete symmetric function
indexed by a list of integers.

\begin{verbatim}
sage: s(st[3,2].character_to_frobenius_image(8))
s[3,3,2]
sage: h(ht[3,2].character_to_frobenius_image(6))
h[3,2,1]
\end{verbatim}

It is this operation that is the origin of our definitions since
in the beginning of our investigations
the $\hht$ and $\st$ basis for us were the pre-images of
the Schur and complete symmetric functions in the $\phi_n$ map.
\end{section}


\begin{thebibliography}{99}

\bibitem[BDO]{BDO}  C. Bowman, M. De Visscher, R. Orellana,  \textit{The partition algebra and the Kronecker coefficients}. Trans. Amer. Math. Soc. 367 (2015), no. 5, 3647Ð3667.

\bibitem[BOR]{BOR} E. Briand, R. Orellana, M. Rosas, 
\textit{The stability of the Kronecker product of Schur functions},
J. Algebra 331 (2011), 11--27.

\bibitem[CG]{CG} C. Chauve, A. Goupil,
\textit{Combinatorial Operators for Kronecker Powers of Representations of $S_n$},
S\'em. Loth. Comb.,
Issue 54, (2005) The Viennot Festschrift, 13 pp.

\bibitem[CF]{ChurchFarb} T. Church, B. Farb,
\textit{Representation theory and homological stability},
Advances in Mathematics,
Volume 245, 1 (2013), Pages 250--314.


\bibitem[GG]{GG} A. Garsia, A. Goupil,
\textit{Character Polynomials, their q-Analogs and the Kronecker Product},
Elec. J. Comb., Volume 16, Issue 2 (2009)
(The Bj\"orner Festschrift volume), \#R19.

\bibitem[Hal]{Halverson} T. Halverson,
\textit{Characters of the Partition Algebras},
J. Algebra,
Volume 238, Issue 2, 15 April 2001, Pages 502--533.

\bibitem[HalRam]{HalversonRam} T. Halverson, A. Ram,
\textit{Partition Algebras},
Europ. J. of Combinatorics 26 (2005) 869--921.

\bibitem[Jon]{Jo} V.F.R.  Jones,  \textit{The  Potts  model  and  the  symmetric  group},  in:  Subfactors:  Proceedings  of  the  Taniguchi Symposium  on  Operator  Algebras  (Kyuzeso,  1993),  River  Edge,  NJ,  World  Sci.  Publishing,  1994, pp. 259Ð267


\bibitem[KT]{KT} K. Koike and I. Terada, 
\textit{Young-diagrammatic methods for the representation theory
of the classical groups of type $B_n$, $C_n$, $D_n$}, 
J. Algebra {\bf 107} (1987) 466--511.

\bibitem[Lascoux]{Lascoux} A. Lascoux,
\textit{Symmetric Functions}, Course about symmetric functions, 
given at Nankai University, October-November 2001, available at
{\tt http://www.emis.de/journals/SLC/wpapers/s68vortrag/ALCoursSf2.pdf}.

\bibitem[Mac]{Mac} I.~G.~Macdonald,
\newblock \textit{Symmetric Functions and Hall Polynomials},
\newblock Second Edition, Oxford University Press,
second edition, 1995.


\bibitem[Ma1]{Ma1}
P. Martin,
\textit{Representations of graph Temperley-Lieb algebras},
Publ. Res. Inst. Math. Sci., 26 (1990), pp. 485--503.

\bibitem[Ma2]{Ma2}
P. Martin,
\textit{Potts Models and Related Problems in Statistical
Mechanics}, World Scientific, Singapore (1991).

\bibitem[Ma3]{Ma3}
P. Martin,
\textit{Temperley-Lieb algebras for nonplanar statistical
mechanics--The partition algebra construction},
J. Knot Theory Ramifications, 3 (1994), pp. 51--82.

\bibitem[Ma4]{Ma4}
P. Martin,
\textit{The structure of the partition algebras},
J. Algebra, 183 (1996), pp. 319--358.

\bibitem[Murg]{Murg} F. D. Murnaghan,
\textit{The characters of the symmetric group},
Amer. J. of Math., 59(4):739--753, 1937.

\bibitem[OZ]{OZ} R. Orellana, M. Zabrocki,
\textit{Symmetric group characters
as symmetric functions}, {\tt arXiv:1605.06672}.

\bibitem[Rosas]{Rosas} M. Rosas,
\textit{A combinatorial overview of the theory of MacMahon
symmetric functions
and a study of the Kronecker product of Schur functions},
Ph.D. Thesis, Brandeis University, 2000.

\bibitem[Sagan]{Sagan} B. Sagan, \textit{The symmetric group}.
Representations, combinatorial algorithms, and
symmetric functions, 2nd edition,  
Graduate Text in Mathematics 203.  Springer-Verlag, 2001.
 xvi+238 pp.

\bibitem[Sp]{Sp} W. Specht, 
\textit{Die Charactere der symmetrischen Gruppe},
Math. Zeitschr. 73, 312--329, 1960.

\bibitem[Stanley]{Stanley} R.~Stanley, 
\textit{Enumerative Combinatorics, Vol. ~2},
Cambridge University Press, 1999.

\bibitem[sage]{sage} W.\thinspace{}A. Stein et~al.
\newblock \textit{{S}age {M}athematics {S}oftware
({V}ersion 4.3.3)},
The Sage Development Team, 2010, {\tt http://www.sagemath.org}.

\bibitem[sage-combinat]{sage-co}
The {S}age-{C}ombinat community.
\newblock  \textit{{{S}age-{C}ombinat}:
enhancing Sage as a toolbox for
computer exploration in algebraic combinatorics},
{{\tt http://combinat.sagemath.org}}, 2008.

\bibitem[Sund]{Sund} S. Sundaram. 
\textit{The Cauchy identity for $Sp(2n)$}. 
J. Combin. Theory Ser. A, 53(2):209--238, 1990.

\bibitem[Weyl]{Weyl}
H. Weyl
\textit{Classical Groups, Their Invariants and Representations},
Princeton Mathematical Series, Princeton Univ. Press,
Princeton (1946).
\end{thebibliography}
\end{document}